\documentclass{article}
\usepackage{graphicx} 
\usepackage{amsmath}
\usepackage{amssymb}
\usepackage{authblk}
\usepackage{orcidlink}
\usepackage{graphicx}
\usepackage{url}

\newtheorem{theorem}{Teorema}[section]

\newtheorem{definition}[theorem]{Definition}
\newtheorem{example}[theorem]{Example}

\title{\textbf{Topological Correlation}}
\author{Isabella Mastroianni\,\orcidlink{0009-0002-9866-3648}}
\author{Ulderico Fugacci\,\orcidlink{0000-0003-3062-997X}}

\affil{\textit{Institute of Applied Mathematics and Information Technologies “Enrico Magenes”, National Research Council, Genova, Italy}}
\date{June 2025}

\begin{document}

\maketitle

\begin{abstract}
    We introduce two novel concepts, \textit{topological difference} and \textit{topological correlation}, that offer a new perspective on the discriminative power of multiparameter persistence. The former quantifies the discrepancy between multiparameter and monoparameter persistence, while the other leverages this gap to measure the interdependence of filtering functions. Our framework sheds light on the expressive advantage of multiparameter over monoparameter persistence and suggests potential applications.
\end{abstract}

\section{Motivation}
Persistent homology is a fundamental tool in topological data analysis, providing insights into the shape of data by tracking the evolution of homological features across a filtration.
While monoparameter persistence has been widely studied and applied \cite{edelsbrunner2010computational, carlsson2021topological, carlsson2009topology, DONUT}, multiparameter persistence has gained increasing attention due to its enhanced discriminative power \cite{carlsson2007theory,kerber2020multi}. However, the reasons behind this greater expressiveness remain an open question, which in turn has contributed to the comparatively lower adoption of multiparameter persistence in actual applications.

Our research seeks to rigorously characterize this more insightful discriminative power by investigating, both from an intuitive and theoretical perspective, how multiparameter persistence encodes richer structural information compared to its single-parameter counterpart.

More specifically, this work represents a first step in this direction by introducing an approach to \textit{quantitative} measure the informative gap through the concept of \textbf{topological difference}. A natural next step is to explore a \textit{qualitative} perspective, aiming at a deeper understanding of its essential meaning. First evidences suggest that this gap may be inherently related to the notion of correlation between filtering functions. As a starting point, we investigate how topological difference can be systematically leveraged to define such a measure of correlation/dependence. We focus on comparing monoparameter and biparameter persistence, for the sake of simplicity and because scalar functions can be effectively compared by examining them in pairs.

To clarify in which sense multiparameter persistence guarantees greater discriminative power compared to monoparameter theory, we present a simple yet representative case.

\begin{example}\label{es:filtrazione S1}
    Let $X$ be the unitary circle centered in the origin and consider the following functions $f,g:X\rightarrow\mathbb{R}, \Phi_1,\Phi_2:X\rightarrow\mathbb{R}^2$: \[ f(x,y):= x \quad 
			g(x,y):= y, \quad 
			\Phi_1:=(f,g), \quad \Phi_2:=(f,f)\]
            
    We obtain the filtrations shown in Figure \ref{fig:filtrazione S1} (a)-(b) and than compute the corresponding persistence homology modules; we will consider the $0$-degree ones, which will be denoted as follows: $M_f, M_g, M_{\Phi_1}$ and $M_{\Phi_2}$. The tables in Figure \ref{fig:filtrazione S1} (c)-(d) display the dimensions of the $0$-degree homology groups\footnote{Recall that they are $\mathbb{K}$-vector spaces when working with coefficients in a field $\mathbb{K}$.} of the sublevel sets associated with each entry.

    Since $M_f\cong M_g$, the monoparametric modules cannot distinguish $\Phi_1$ from $\Phi_2$, whereas $M_{\Phi_1}\ncong M_{\Phi_2}$,\footnote{For instance, Figure \ref{fig:filtrazione S1} (c)-(d) shows that at $(-1,-1)$ the dimensions of $M_{\Phi_1}$ and $M_{\Phi_2}$ are $0$ and $1$ respectively.} making a distinction possible.
 
    \begin{figure}[!hbtp]
        \centering
        \includegraphics[width=\textwidth]{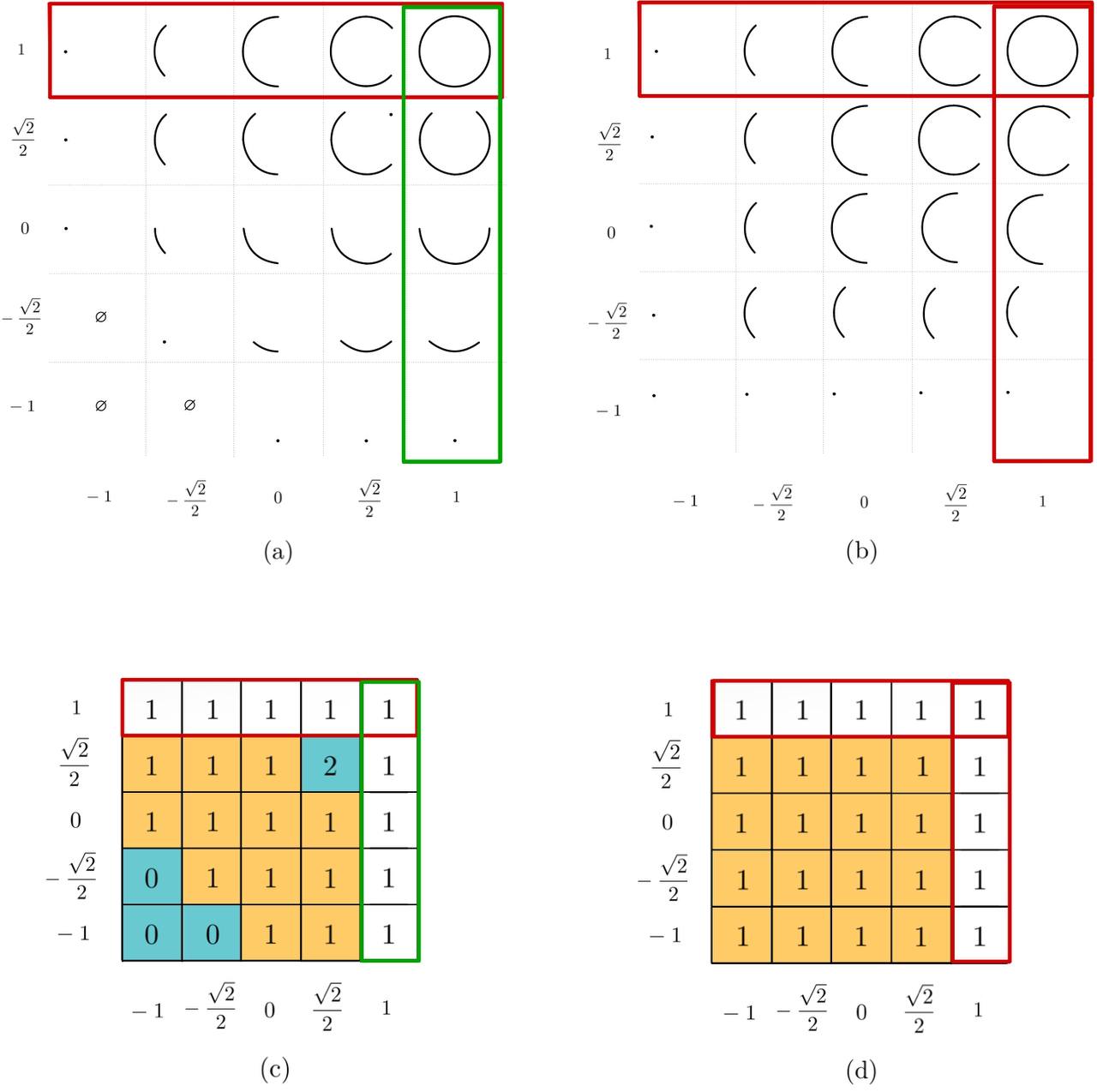}
        \caption{Filtrations induced by $\Phi_1$ (a), $\Phi_2$ (b). Dimensions of $M_{\Phi_1}$ (c), $M_{\Phi_2}$ (d). Filtrations induced by $f$, $g$ and dimensions of $M_f$, $M_g$ are framed respectively in red and green.}
        \label{fig:filtrazione S1}
    \end{figure}

\end{example}

Similar constructions can be easily devised, even with the same level of simplicity.

\section{Definitions}

Throughout the rest of this discussion, we will work under the following hypoteses (as in \cite{frosini2023matching}):
\begin{itemize}
    \item $X$ a finitely triangulable space;
    \item all the considered functions will be continuous;
    \item for every $k$ integer and $u\in\mathbb{R}^d$ ($d=1$ or $d=2$), the \v{C}ech homology group $H_k(X_u)$ of the sublevel set of $X$ with respect to any function will be finitely generated.
\end{itemize}

We introduce the definition of topological difference to quantify the informative gap between a multiparametric and a monoparametric approach to function comparison.

\begin{definition}
    Given two functions $ \Phi_1:=(f_1,g_1),\Phi_2:=(f_2,g_2):X\rightarrow \mathbb{R}^2 $, we define the \textbf{topological difference} between $\Phi_1, \Phi_2$ as \[\Delta(\Phi_1,\Phi_2):= d_{match}(\Phi_1,\Phi_2)-\max\{d_B(f_1,f_2),d_B(g_1,g_2)\},\]
    where $d_B$ is the bottleneck distance and $d_{match}$ is the matching distance that incorporates also horizontal and vertical filtering lines \cite{BIASOTTI20111735,frosini2023matching}.
\end{definition}

Observe that this quantity is non-negative, as the matching distance is defined as the supremum of the bottleneck distances across a foliation of the considered biparametric filtrations, which also accounts for horizontal and vertical slopes corresponding to the original filtrations induced by the components.

Building on this concept, we also introduce a measure of correlation between the components of a biparametric function $\Phi:=(f,g):X\rightarrow \mathbb{R}^2 $. The key idea is to consider $F:=(f,f)$ and $G:=(g,g)$, compute the largest possible discrepancy \, ---\, namely the maximum between $\Delta(\Phi,F)$ and $\Delta(\Phi,G)$ \, ---\, and subsequently apply a suitable normalization. This leads us to the following definition.

\begin{definition}
    Given $\Phi:=(f,g):X\rightarrow \mathbb{R}^2 $ we define the \textbf{topological correlation} $\Delta(\Phi)$ between its components as
    \[ \Delta(\Phi):=
		\begin{cases}
			\frac{2\max\{\Delta(\Phi,F),\Delta(\Phi,G)\}}{\Delta(\Phi,F)+\Delta(\Phi,G)}-1 & if\ \Delta(\Phi,F)+\Delta(\Phi,G) \neq 0, \\
			1 & otherwise.
		\end{cases}\]
\end{definition}

Notice that the condition $\Delta(\Phi,F)+\Delta(\Phi,G) = 0$ is equivalent to $\Delta(\Phi,F)=0=\Delta(\Phi,G)$, since topological difference is non negative. This occurs when there is no discrepancy between the information provided by the monoparameter filtrations and the biparameter filtration.

It is not so difficult to check that, as expected, $\Delta(\Phi)$ takes values in $[0,1]$ and that, in accordance with intuition, the following properties are satisfied:

\begin{enumerate}
    \item If $\Phi:=(f,f)$, then $\Delta(\Phi)=1$;
    \item $\Delta(\Phi)=0$ if and only if $\Delta(\Phi,F)=\Delta(\Phi,G)\neq 0$;
    \item $\Delta(\Phi)=1$ if and only if $d_{match}(\Phi,F)=d_B(f,g)=d_{match}(\Phi,G)$.
\end{enumerate}

Since we aim at adopting these tools to possible applications, we put emphasis on scalar function comparison, which will be illustrated in the next section.

\section{Applications}

The notions just introduced provide the basis for formulating a measure of dependence between functions defined on different domains.

Let $X_1,...,X_N$ be some topological spaces and $\mathcal{F}:=\{f_i\},\mathcal{G}:=\{g_i\}$ two collection of $\mathbb{R}$-valued functions defined on $X_i$, all satisfying the hypotheses stated in the previous section.

We can now consider $\Phi_i := (f_i, g_i): X_i \rightarrow \mathbb{R}^2$ in order to obtain a collection of values $\{\Delta(\Phi_i)\}_{i=1,...,N}$ so that, by considering their average, a definition of correlation between $f$ and $g$ is given.

\begin{definition}
    We define the \textbf{topological correlation} between $f$ and $g$ as the arithmetic mean of $\{\Delta(\Phi_i)\}$, namely:  
\[
\operatorname{corr}(\mathcal{F},\mathcal{G}) := \frac{1}{N} \sum_{i=1}^{N} \Delta(\Phi_i).
\]
\end{definition}

The motivation for this setting is due to possible applications. For example, we can consider as the spaces $X_i$ some molecular surfaces and as $f_i$, $g_i$ some functions describing physicochemical property of $X_i$ \cite{raffo2021shrec}. We conclude with the following example.

\begin{example}
    Let us consider Example \ref{es:filtrazione S1}. In such a case, we obtain $d_B(f,g)=0$ which implies that $\Delta(\Phi_1,F_1)=d_{match}(\Phi_1,F_1)$ and $\Delta(\Phi_1,G_1)=d_{match}(\Phi_1,G_1)$. Moreover, by symmetry, we can verify that $ d_{match}(\Phi_1,F)=d_{match}(\Phi_1,G)\geq \dfrac{1}{2}-\dfrac{\sqrt{2}}{4}$. It follows that $\Delta(\Phi_1,F_1)=\Delta(\Phi_1,G_1)\neq 0$  and hence $\Delta(\Phi_1)=0$.
    
    A similar argument applies when considering the projections onto the $x-$axis and $y-$axis, defined on the unit sphere and a torus $S^1\times S^1$.
    
    Therefore, the correlation between the two collection of functions $\mathcal{F}$ and $\mathcal{G}$ given by the $x-$axis and $y-$axis projections in this case is zero.
\end{example}

\section{Advantages, limitation and future perspective}

The theoretical foundations of these concepts are supported by rigorous properties and illustrative examples. However, we acknowledge certain limitations that need to be addressed and refined.

The primary challenge in computing topological difference and, consequently, topological correlation lies in the necessity of first computing the matching distance (or an approximation thereof), which can be highly complex. Currently, we are following the approach described in \cite{BIASOTTI20111735, cerri2011new}; another possible reference is undoubtedly \cite{kerber2020exact}.

Moreover, some aspects of the classical notion of correlation remain unclear in this new setting. For instance, the relevance of a notion of negative correlation, the expected behavior when comparing a function $f$ with its scaled version $k \cdot f$, and the treatment of potential nonlinear correlations, such as those between $x$ and $x^2$, all require further investigation. However, these challenges may not constitute a fundamental drawback when considering our primary goal: determining whether a monoparameter or a multiparameter approach is more suitable in a given context.

Additionally, our current definitions rely on bottleneck and matching distances, which are not true metrics but rather pseudo-metrics. A possible alternative could be the interleaving distance, yet its computation is so demanding that it may not be a practical option.

Finally, at this stage, we have not extensively conducted concrete experiments yet. This is partly because we are working to formalize an intuitive idea that currently lacks a precise definition, making it difficult to establish a basis for comparison. We are nonetheless working to lay the groundwork for such experiments and plan to carry them out in the future.

For all these reasons, the next steps in our research will be the following.

\begin{itemize}
\item Developing and implementing an algorithm for computing the matching distance \, ---\, an aspect we are already actively working on. The current code is available in the following GitHub repository: \href{https://github.com/IsHubolla/Matching-distance}{\textit{Matching-distance}}.
\item Testing our definitions in different contexts, assessing their robustness, and refining them as needed.
\item Exploring equivalent definitions directly in terms of persistence modules \cite{oudot2017persistence}. In particular, given a function $\Phi:=(f,g):X\rightarrow\mathbb{R}^2$, we have the corresponding persistence modules:
\begin{center}
\begin{tabular}{lcr}
$M_f$ & $M_g$ & $M_{\Phi}$.
\end{tabular}
\end{center}
If an appropriate product operation $\star$ on persistence modules could be defined, one could investigate whether the relation $$M_f\star M_g=M_{\Phi}$$ holds, providing a structural interpretation of topological correlation. A preliminary investigation of a such possibility can be found in \cite{mastroianni2025retrievingbiparameterpersistencemodules}.
\end{itemize}

\section*{Acknowledgments}
This work has been developed within the CNR research activities STRIVE DIT.AD022.207.009 and  DIT.AD021.080.001 and within the framework of the projects ``RAISE - Robotics and AI for Socio-economic Empowerment'' - Spoke number 2 (Smart Devices and Technologies for Personal and Remote Healthcare) and 
``NBFC - National Biodiversity Future Center'' - Spoke number 4 (Ecosystem Functions, Services and Solutions).

The authors wish to thank Silvia Biasotti and Marco Guerra for the helpful discussions and suggestions.

\bibliographystyle{plainurl}
\bibliography{refs}

\begin{thebibliography}{10}

\bibitem{BIASOTTI20111735}
Silvia Biasotti, Andrea Cerri, Patrizio Frosini, and Daniela Giorgi.
\newblock A new algorithm for computing the 2-dimensional matching distance between size functions.
\newblock {\em Pattern Recognition Letters}, 32(14):1735--1746, 2011.

\bibitem{carlsson2009topology}
Gunnar Carlsson.
\newblock Topology and data.
\newblock {\em Bulletin of the American Mathematical Society}, 46(2):255--308, 2009.

\bibitem{carlsson2021topological}
Gunnar Carlsson and Mikael Vejdemo-Johansson.
\newblock {\em Topological data analysis with applications}.
\newblock Cambridge University Press, 2021.

\bibitem{carlsson2007theory}
Gunnar Carlsson and Afra Zomorodian.
\newblock The theory of multidimensional persistence.
\newblock In {\em Proceedings of the twenty-third annual symposium on Computational geometry}, pages 184--193, 2007.

\bibitem{cerri2011new}
Andrea Cerri and Patrizio Frosini.
\newblock A new approximation algorithm for the matching distance in multidimensional persistence.
\newblock 2011.

\bibitem{edelsbrunner2010computational}
Herbert Edelsbrunner and John Harer.
\newblock {\em Computational topology: an introduction}.
\newblock American Mathematical Soc., 2010.

\bibitem{frosini2023matching}
Patrizio Frosini, Eloy~M{\'o}sig Garc{\'\i}a, Nicola Quercioli, and Francesca Tombari.
\newblock Matching distance via the extended pareto grid.
\newblock {\em arXiv preprint arXiv:2312.04201}, 2023.

\bibitem{DONUT}
Barbara Giunti, J{\=a}nis Lazovskis, and Bastian Rieck.
\newblock {DONUT}: {D}atabase of {O}riginal \& {N}on-{T}heoretical {U}ses of {T}opology, 2022.
\newblock \url{https://donut.topology.rocks}.

\bibitem{kerber2020multi}
Michael Kerber.
\newblock Multi-parameter persistent homology is practical.
\newblock In {\em TDA \& Beyond}, 2020.

\bibitem{kerber2020exact}
Michael Kerber, Michael Lesnick, and Steve Oudot.
\newblock Exact computation of the matching distance on 2-parameter persistence modules.
\newblock {\em Journal of Computational Geometry}, 11(2):4--25, 2020.

\bibitem{mastroianni2025retrievingbiparameterpersistencemodules}
Isabella Mastroianni, Marco Guerra, Ulderico Fugacci, and Emanuela~De Negri.
\newblock Retrieving biparameter persistence modules from monoparameter ones: a characterization of hook-decomposable persistence modules, 2025.
\newblock URL: \url{https://arxiv.org/abs/2506.14678}, \href {https://arxiv.org/abs/2506.14678} {\path{arXiv:2506.14678}}.

\bibitem{oudot2017persistence}
Steve~Y Oudot.
\newblock {\em Persistence theory: from quiver representations to data analysis}, volume 209.
\newblock American Mathematical Soc., 2017.

\bibitem{raffo2021shrec}
Andrea Raffo, Ulderico Fugacci, Silvia Biasotti, Walter Rocchia, Yonghuai Liu, Ekpo Otu, Reyer Zwiggelaar, David Hunter, Evangelia~I Zacharaki, Eleftheria Psatha, et~al.
\newblock Shrec 2021: Retrieval and classification of protein surfaces equipped with physical and chemical properties.
\newblock {\em Computers \& Graphics}, 99:1--21, 2021.

\end{thebibliography}

\end{document}